\documentclass{svmult}
\font\lou=msbm10 scaled 1000
\def\N{\hbox{\lou N}}

\def\ds{\displaystyle}

\usepackage{graphicx}

\usepackage{amssymb}
\usepackage{dsfont}

\newenvironment{enumerate-roman}{\begin{enumerate}}{\end{enumerate}}

\usepackage{makeidx}         

\usepackage{graphicx}        

\usepackage{multicol}        

\usepackage[bottom]{footmisc}

\makeindex             

\newtheorem{thm}{Theorem}[section]
\newtheorem{lem}[thm]{Lemma}
\newtheorem{prop}[thm]{Proposition}

\newtheorem{defi}[thm]{Definition}

\newcommand{\mod}{{\rm \  mod \ }}

\begin{document}

\baselineskip15pt

\title*{Random Periodic Solutions of Random Dynamical Systems 
}
\titlerunning{Random periodic solutions}
\author{Huaizhong Zhao$^{1}$, Zuo-Huan Zheng$^{1,2}$ }
\authorrunning{H.Z. Zhao and Z.H. Zheng}
\institute{Department of Mathematical Sciences, Loughborough
University, LE11 3TU, UK. \texttt{H.Zhao@lboro.ac.uk}
\and
Institute of Applied
Mathematics, Academy of Mathematics and Systems Sciences, Chinese
Academy of Sciences, Beijing 100080, P. R. China.
\texttt{zhzheng@amt.ac.cn}}
\maketitle

\begin{abstract}
In this paper, we give the definition of the random periodic solutions
of random dynamical systems. We prove the existence of such periodic solutions for a $C^1$ perfect cocycle
on a cylinder using a random invariant set, the Lyapunov exponents
and the pullback of the cocycle.
\end{abstract}
\vskip0.4cm



\section{Introduction}

Similar to the deterministic dynamical systems, in stochastic dynamical systems, the problem
of the long time or infinite horizon behaviour is a fundamental problem
occupying a central place of research. Studies of dynamic properties of such systems, both
deterministic and stochastic, usually involve an appropriate definition of a steady state
(viewed as a dynamic equilibrium) and conditions that guarantee its existence and
local or global stability.

Fixed points or periodic solutions capture the intuitive
idea of a stationary state or an equilibrium of a dynamical system.
For a deterministic dynamical system $\Phi_t: {\cal H}\to {\cal H}$ over time $t\in I$, where ${\cal H}$ is
the state space, $I$ is the set of all real numbers, or discrete real numbers, a fixed point
is a point $y\in {\cal H}$ such that
\begin{eqnarray}\label{fp1}
\Phi_t(y)=y {\rm \ for \ all} \  t\in I,
\end{eqnarray}
and a periodic solution is a periodic function $\psi: I\to {\cal H}$ with period $T\neq 0$ such that
\begin{eqnarray}\label{pe1}
\psi(t+T)=\psi (t) \ {\rm and }\ \Phi_t(\psi(t_0))=\psi(t+t_0)
{\rm \ for \ all} \  t, t_0\in I.
\end{eqnarray}
However, in many real problems, such kinds of the deterministic
equilibria actually don't exist due to the presence of random
factors. For stochastic dynamical systems, it would be reasonable to
say that stationary states are not actually steady states in the
sense of (\ref{fp1}) and periodic solutions in the sense of
(\ref{pe1}). Due to the fact that the random external force pumps
to the system constantly, the relation (\ref{fp1}) or (\ref{pe1})
breaks down.

Let $\Phi: \Omega\times I\times {\cal H}\to {\cal H}$ be a measurable random
dynamical system on a measurable space $({\cal H},{\cal B})$ over a metric
dynamical system $(\Omega, {\cal F}, P,(\theta_t)_{t\geq 0})$. Then
a stationary solution of $\Phi$ is a ${\cal F}$ measurable random
variable $Y:\Omega \to {\cal H}$ such that (c.f. Arnold \cite{arnold})
\begin{eqnarray} \Phi(\omega,t,Y(\omega))=Y(\theta _t\omega){ \ \rm
for \ all \ } t\in I\  a.s.
\end{eqnarray}
The concept of the stationary random point of a random dynamical
system is a natural extension of the equilibrium or fixed point in
deterministic systems. Such a random fixed point for SPDEs is
especially interesting. It consists of  infinitely many random
moving invariant surfaces on the configuration space. It is a more
realistic model than many deterministic models as it demonstrates
some complicated phenomena such as turbulence. Finding such
stationary solutions for SPDEs is one of the basic problems. But the
study of random cases is much more difficult and subtle, in contrast
to deterministic problems. This ``one-force, one-solution" setting
describes the pathwise invariance of the stationary solution over
time along $\theta$ and the pathwise limit of the random dynamical
system. Their existence and/or stability for various stochastic
partial differential equations have been under active study recently
(Duan, Lu and Schmalfuss \cite{duan}, E, Khanin, Mazel and Sinai
\cite{sinai2}, Mohammed, Zhang and Zhao \cite{mohammed},
Zhang and Zhao \cite{zhang}).

Needless to say that the random periodic solution is another
fundamental concept in the theory of random dynamical systems.
According to our knowledge, such a notion did not exist in
literatures. In this paper, we will carefully define such a notion
and give a sufficient condition for the existence on a cylinder
$S^1\times R^{d}$. To see the motivation for such a definition,
let's first note two obvious but fundamental truth in the
definition of periodic solution (\ref{pe1}) of the deterministic
systems when $I$ is the set of real numbers:-

(i) The function $\psi$ (given in the parametric form here) is a closed curve in
the phase space;

(ii) If the dynamical system starts at a point on the closed curve, the
orbit will remain on the same closed curve.

\noindent But note, in the case of stochastic dynamical systems,
although the function $\psi$ may still be a periodic function, one would
expect that $\psi$ depend on $\omega$. In other words, we would expect infinitely
many periodic functions $\psi^{\omega}$, $\omega\in \Omega$.
Moreover, even the random
dynamical system starts at a point on the curve $\psi^{\omega}$, it will
not stay in the same periodic  curve when time is running. In fact, the periodic
curve actually is not the orbit of the random dynamical system, but
the random dynamical system will move from one periodic curve
$\psi^{\omega}$ to another periodic curve $\psi^{\theta_t\omega}$ at time
$t\in I$. Now we give the following definition:

\begin{defi} \label{revision2} A random
periodic solution is an ${\cal F}$-measurable periodic function
$\psi:\Omega\times I\to {\cal H}$ of period $T$ such that
\begin{eqnarray}\label{pe2} \psi^{\omega}(t+T)=\psi ^{\omega}(t) \ {\rm and }\
\Phi_t^{\omega}(\psi^{\omega}(t_0))=\psi^{\theta _t \omega}(t+t_0) {\rm \ for
\ all} \ t, t_0\in I.
\end{eqnarray}
\end{defi}

Comparing to the stationary solution, which consists of infinitely
many single points on the phase space, the random periodic
solution consists of infinite many random moving periodic curves on
the phase space. It describes more complex random nonlinear
phenomenon than a stationary solution. The closed curve on the phase
space is a pathwise invariant set over the time and the
$P$-preserving transformation $\theta: I\times \Omega\to \Omega$
along the orbit defined by the random dynamical system. To
understand and give the existence of such a solution is an
interesting problem. The systematic study of the deterministic
periodic solutions of the deterministic dynamical systems began in
Poincar\'e in his seminal work \cite{poincare}. The
Poincare-Bendixson Theorem has been very useful in the study of the
periodic solutions (\cite{bendixson}). 
Periodic solutions have been studied for many important systems arising in 
numerous physical problems e.g. van der Pol equations (van der Pol \cite{vanderpol}),
Li\'enard equations (Li\'enard \cite{lienard}, Filipov \cite{fili}, Zheng 
\cite{zheng3}).
 Now, after over a century,
this topics is still one of the most interesting nonlinear phenomena
to study in the theory of deterministic dynamical systems. The
periodic solutions have occupied a central place in dynamical
systems such as in the study of the bifurcation theory (V. Arnold
\cite{varnold}, Chow and Hale \cite{chow}, Andronov \cite{Andronov}, Li and Yorke \cite{li} to
name but a few). The Hilbert's 16th problem involves determining the
number and location of limits cycles for autonomous planar
polynomial vector field. This problem still remains unsolved
although  much progress has been made. See Ilyashenko
\cite{ilyashenko} for a summary. Needless to say that the study of
the random periodic solution is more difficult and subtle.
The extra essential difficulty comes from the fact that the
trajectory of the random dynamical systems starting at a point on
the periodic curve does not follow the periodic curve, but moves
from one periodic curve to another one corresponding to different
$\omega$. If one looks at a family of trajectories (therefore forms a
map or a flow) starting from different points in the closed curve
$\psi^{\omega}$, then the whole family of trajectories at time $t\in
I$ will lie on a closed curve corresponding to $\theta _t\omega$.
This is essentially different from the deterministic case. Of
course, the definition (\ref{pe1}) in the deterministic case is a
special case of the definition (\ref{pe2}) of the random case.

In section 2, we will give the definition of
the random periodic  solutions of a random dynamical system on a cylinder. In
section 3, we will prove the existence of the finite number of periodic solution
assuming a contraction condition using Lyapunov exponent
near the attractor.  We use the compactness
argument and obtain a finite number of open covering and eventually to
prove the existence of the finite number of closed curve. But this does not give
the estimate on how many periodic solutions and the winding numbers. They are
interesting problems to investigate in the future. Moreover, we will prove that the winding number of 
each closed orbit, the period of each periodic solution
and the number of periodic solutions are invariant under the perturbation of noise.

\section {The notion of periodic invariant  solutions on cylinders and an example}

The extension of the notion  of a periodic solution on a cylinder to the random case
is given as follows (see Fig.1), where $I$ is either $[0,\infty)$, or $(-\infty,0]$, or
$(-\infty,+\infty)$. Note here there is a natural parameter $s\in S^1$ for a closed curve on
$S^{1}\times R^{d}$. Some cases on $R^{d+1}$ can be transformed to the cylinderic case. In the next section,
we will give the existence of the periodic solutions of the random dynamical systems on a cylinder.
The following definition gives more information about the winding number of
 the closed curves. Except for the winding number, on a cylinder, Definition \ref{revision1}
 is equivalent to Definition \ref{revision2} if we reparameterize the natural parameter $s$ using time
 by putting $s=s(t)$. But on a cylinder,  it is more convenient to use the natural 
 parameter $s$. 

\begin{defi}\label {revision1}
Let $\varphi^{\omega}: R\to R^{d}$ be a continuous periodic
function of period $\tau\in \N$ for each $\omega\in \Omega$. Define
$L^{\omega}={\rm graph}(\varphi^{\omega})=\{(s \mod
1,\varphi^{\omega}(s)):s\in R^1\}$. If $L^{\omega}$ is invariant
with respect to the random dynamical system
$\Phi:\Omega\times {I}\times S^1\times R^{d}\to
S^1\times R^{d}$, i.e.
$\Phi^{\omega}(t)L^{\omega}=L^{\theta_t\omega}$, and  there exists a
minimum $T>0$ (or maximum $T<0$) such that for any $s\in [0,\tau)$,
\begin{eqnarray}
\Phi^{\omega}(T,(s\mod 1,\varphi^{\omega}(s)))
=(s\mod 1,\varphi^{\theta_{T}\omega}(s)),
\end{eqnarray}
for almost all $\omega$, then it is said that $\Phi$ has a random
periodic solution of period $T$ with random periodic curve $L^{\omega}$ 
of winding number $\tau$. \end{defi} 

  It is easy to see that
$\delta _{L^{\omega}}(dx)P(d\omega)$ is an invariant measure of the
skew-product $(\theta, \Phi): I\times \Omega\times (S^1\times R^d)\to \Omega \times (S^1\times R^d)$.
Needless to say, an invariant measure may not give a random periodic solution. The
support of an invariant measure may be a random fixed point (stationary solution), or
random periodic solutions, or a more complicated set.

The periodic invariant  orbit is a new concept in the literature. We believe it has some
importance in random dynamical systems, for example, it can be
studied systematically to establish the Hopf bifurcation theory of
stochastic dynamical systems. This is not the objective of this
paper, we will study this problem in future publications. But here
in order to illustrate the concept, as a simple example, we
consider the random dynamical system generated by a perturbation to
the following deterministic ordinary differential equation in $R^2$:
\begin{eqnarray}\
\cases {{{dx(t)}\over {dt}}=x(t)-y(t)-x(t)(x^2(t)+y^2(t)),\cr
{{dy(t)}\over {dt}}=x(t)+y(t)-y(t)(x^2(t)+y^2(t)).}
\end{eqnarray}
It is well-known that above equation has a limit cycle
\begin{eqnarray*}
x^2(t)+y^2(t)=1.
\end{eqnarray*}
\vskip20pt

\begin{figure}[t]
   \vspace{0.5in}
\includegraphics[width=4in]{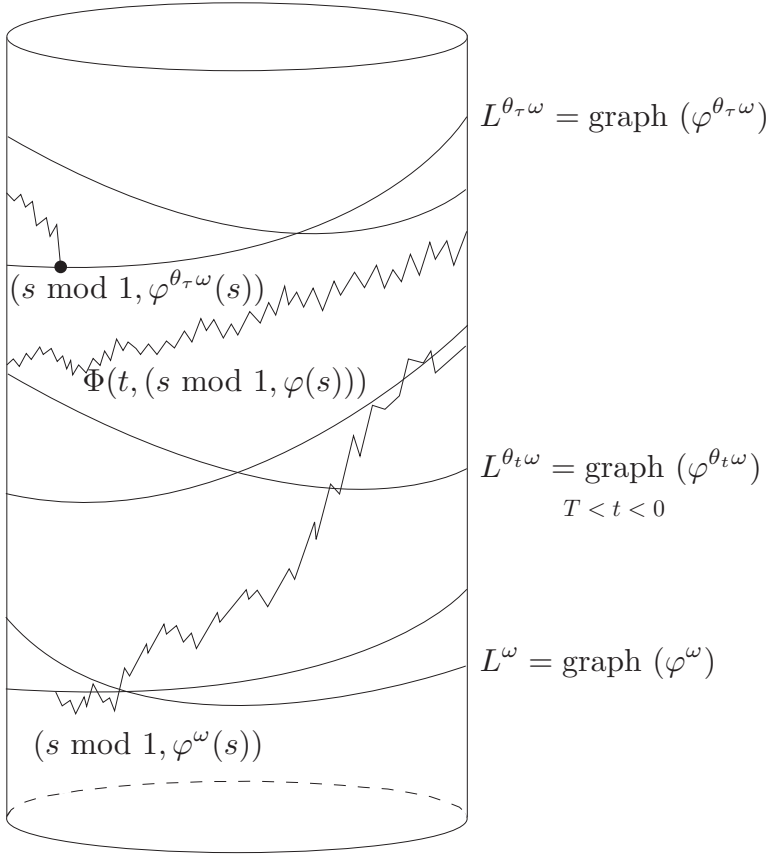}
   \caption{Random periodic orbit of period $T$ and
winding number $\tau=2$.}
\end{figure}

Consider a random perturbation
\begin{eqnarray}\label{sode1}
\cases {dx=(x-y-x(x^2+y^2))dt+x\circ dW(t),\cr
dy=(x+y-y(x^2+y^2))dt+y\circ dW(t).}
\end{eqnarray}
Here $W(t)$ is a one-dimensional motion on the canonical probability space $(\Omega,{\cal F},P)$
with the $P$-preserving map $\theta$ being taken to the shift operator
$(\theta_t\omega)(s)=W(t+s)-W(t)$.
Using polar coordinates
\begin{eqnarray*}
x=\rho \cos2\pi \alpha, \ \ y=\rho \sin2\pi\alpha,
\end{eqnarray*}
then we can transform Eq. (\ref{sode1}) on $R^2$ to the following equation on the cylinder $[0,1]\times R^1$: 
\begin{eqnarray}\label{peq}
\cases {d\rho(t)=(\rho(t)-\rho^3(t))dt+\rho(t)\circ dW(t),\cr
d\alpha={1\over 2\pi}dt.}
\end{eqnarray}
This equation has a unique close form solution as follows:
\begin{eqnarray*}
\rho(t, \alpha _0,\rho_0,\omega)={\rho_0e^{t+W_t(\omega)}\over
(1+2\rho_0^2\int_0^te^{2(s+W_s(\omega))}ds)^{1\over2}},\ \
\alpha (t,\alpha _0,\rho_0,\omega)=\alpha_0+{t\over 2\pi}.
\end{eqnarray*}
It is easy to check that
\begin{eqnarray*}
\rho^*(\omega)=(2\int_{-\infty}^0e^{2s+2W_s(\omega)}ds)^{-{1\over2}}
\end{eqnarray*}
is the stationary solution of the first equation of (\ref{peq}) i.e.
\begin{eqnarray*}
\rho(t,\alpha _0,\rho^*(\omega),\omega)=\rho^*(\theta_t\omega)
\end{eqnarray*}
and
\begin{eqnarray*}
\Phi(t,\omega)(\alpha _0,\rho_0)=(\alpha_0+{t\over 2\pi}\ mod\
1,\rho(t,\alpha _0,\rho_0,\omega))
\end{eqnarray*}
defines a random dynamical system on the cylinder $[0,1]\times R^1$:
$\Phi(t,\omega)=(\Phi_1(t,\omega),\Phi_2(t,\omega)):
[0,1]\times
R^1\longrightarrow [0,1]\times R^1$.
Define
\begin{eqnarray*}
L^\omega=\{(\alpha,\rho^*(\omega)):\ 0\leq\alpha \leq 1\}.
\end{eqnarray*}
Then
\begin{eqnarray*}
L^{\theta_t\omega}=\{(\alpha ,\rho^*(\theta_t\omega)):\ 0\leq\alpha
\leq 1\}.
\end{eqnarray*}
It is noticed that
\begin{eqnarray*}
\Phi(t,\omega)L^\omega &=&\{(\alpha+{t\over 2\pi}\ mod \
1,\rho^*(\theta_t\omega)):\ 0\leq\alpha \leq 1\}\\
&=&\{(\alpha ,\rho^*(\theta_t\omega)):\ 0\leq\alpha \leq 1\}.
\end{eqnarray*}
Therefore
\begin{eqnarray*}
\Phi(t,\omega)L^\omega=L^{\theta_t\omega},
\end{eqnarray*}
i.e. $L^\cdot$ is invariant under $\Phi$. Moreover
\begin{eqnarray*}
\Phi(2\pi,\omega)(\alpha ,\rho^*(\omega))=(\alpha ,\rho^*(\theta_{2\pi}\omega)).
\end{eqnarray*}
Therefore the random dynamical system has a random solution of 
period $2\pi$ with invariant closed 
curve $L^{\omega}$ of winding number 1 on 
$[0,1]\times R^1$. Now we can transform the periodic solution back to $R^2$. For this define
 for $(x,y)\in R^2,\ x=\rho \cos2\pi\alpha,\ y=\rho \sin2\pi\alpha$
\begin{eqnarray*}
&&
\tilde{\Phi}(t,\omega)(x,y)\\
&=&(\Phi_2(t,\omega)(\alpha,\rho)\cos(2\pi \Phi_1(t,\omega)(\alpha,\rho)),\
\Phi_2(t,\omega)(\alpha,\rho)\sin(2\pi \Phi_1(t,\omega)(\alpha,\rho))),
\end{eqnarray*}
and
\begin{eqnarray*}
\psi^{\omega}(t) =(\rho^*(\omega)\cos(2\pi \alpha +t),\rho^*(\omega)\sin
(2\pi \alpha +t)).
\end{eqnarray*}
It is obvious that
\begin{eqnarray*}
\psi^{\omega}(2\pi+t)=\psi^{\omega}(t),
\end{eqnarray*}
and
\begin{eqnarray*}
\tilde \Phi(t,\omega)\psi^{\omega}(0)&=&\tilde \Phi(t,\omega)(\rho^*(\omega)\cos(2\pi \alpha),\rho^*(\omega)\sin
(2\pi \alpha))\\
&=& (\rho^*(\theta _t\omega)\cos(2\pi \alpha +t),\rho^*(\theta _t\omega)\sin
(2\pi \alpha +t))\\
&=&\psi^{\theta _t\omega}(t).
\end{eqnarray*}
From this we can tell that the random dynamical system generated by
the stochastic differential equation (\ref{sode1}) has a random periodic
solution $\psi^{\omega}: \Omega \times I\to R^2$ defined above. Moreover if $x(0)^2+y(0)^2\neq 0$, then
$$
x^2(t,\theta(-t,\omega))+y^2(t,\theta(-t,\omega))\to
\rho^*(\omega)^2$$ as $t\to \infty$.

\section{The existence of random periodic  solutions}

Consider a continuous time differentiable random dynamical system $\gamma$
over a metric dynamical system $(\Omega, {\cal F}, P,\theta)$ on a
cylinder $S^1\times R^{d}$. Though the precise definition given below is standard, it is crucial to the development of this article (c.f. \cite{arnold}, \cite{mohammed}). 

\begin{defi} A $C^1$ perfect cocycle is a ${\cal B}((-\infty,+\infty))\otimes {\cal F}\otimes {\cal
B}(S^1\times R^{d}), {\cal B}(S^1\times R^{d})$ measurable random field  $\gamma:
(-\infty,+\infty)\times \Omega \times S^1\times R^{d}\to S^1\times
R^{d}$ satisfying the following conditions:

 (a) for each $\omega\in
\Omega$, $\gamma (0,\omega)=Id$,

 (b) for each $\omega\in
\Omega$, $\gamma _z^{\omega}(t_1+t_2)=\gamma _{\gamma
_z^{\omega}(t_1)}^{\theta _{t_1}\omega}(t_2)$ for all $z\in
S^1\times R^{d}$ and $t_1,t_2\in R$, 

(c) for each $\omega\in \Omega$, the mapping $\gamma
_{\cdot}^{\omega}(\cdot): (-\infty,+\infty)\times S^1\times R^{d}\to
S^1\times R^{d}$ is continuous,


(d) for each $(t,\omega)\in (-\infty,+\infty)\times \Omega$, the mapping
$\gamma _{\cdot}^{\omega}(t): S^1\times R^{d}\to S^1\times
R^{d}$ is a $C^1$ diffeomorphism.

\end{defi}

In the following, we will assume there is an invariant set $X^{\omega}\subset S^1\times Y^{\omega}$,
for some compact set $Y^{\omega}$, and 
 use a random quasi-periodic winding system on $S^1\times Y^{\omega}$
  to characterize the invariant set $X^{\omega}$ 
further and to study the existence of random periodic curves on the cylinder. 
We will prove under the following three conditions, the invariant set $X^{\omega}$ consists of 
a finite number of periodic curves. The Lyapunov 
exponent and the pullback are key techniques in our approach. So
it is essential to assume the random dynamical systems are $C^1$ perfect cocycles.
It is known from the works of many mathematicians that cocycles can be generated from some random differential equations and stochastic 
differential equations (c.f. Arnold \cite{arnold}).
\vskip10pt

{\bf Condition (i)} {\it Assume that there exists a random compact
subset Y of $R^{d}$ and $t_1>0$ such that for any $s \in S^1$
and $y=(x_1,x_2,\cdots,x_d)\in Y^{\omega}$ (denote
$z_0=(s,x_1,\cdots,x_d))$, 
$\gamma_{z_0}^{\omega}(t_1)=(s',x_1',\cdots,x_d')$ with
\begin{eqnarray}\label{march18a}
s'=s + 1,
\end{eqnarray}
and
\begin{eqnarray}\label{march18b}
y'=(x_1',x_2',\cdots,x_d')\in Y^{\theta _{t_1}\omega}.
\end{eqnarray}
We assume there is an invariant compact set $X^{\omega}\subset S^1\times
Y^{\omega}$, that is to say $
X^{\omega}$ satisfies
\begin{eqnarray}\label{is1}
\gamma ^{\omega}(t)X^{\omega}=X^{\theta_t \omega},
\end{eqnarray}
 for all $t\in R$ and the
projection of $X^{\omega}$ to the subspace $S^1$ is $S^1$, there exist a constant $b^*>0$ such that 
$|X^{\omega}|<b^*$ where $|X^{\omega}|$ denotes the diameter of $X^{\omega}$ in the y-direction.
}
\vskip5pt

Note that
$t_1$ is the time that a particle on $S^1$ rotate a full circle, no matter where it starts, 
(\ref{march18b}) is just an invariant condition about the random compact subset $Y^{\omega}\subset R^d$.

In the following, we denote $\hat \theta\omega=\theta _{t_1}\omega$.
Under the above assumption, a random winding system can be defined from
the continuous time $C^1$ random dynamical system $\gamma$ by
letting
$$
h(s)=s^{\prime}$$ and
$$
g^{\omega}(s,y)=y^{\prime}.
$$
 Consider the following discrete time random winding
system
\begin{eqnarray}
\label{(10)}
 \left \{ \begin{array}{l}
 s_{n+1} = h(s_n)  {\mod} 1 , s_n \in S^1\cr
 y_{n+1} = g ^{\omega}(s_n,y_n),\hspace{0.45cm} y_n \in R^{d},
 \end{array}
 \right.
\end{eqnarray}
where $S^1$ is the unit circle, $g:\Omega\times S^1\times R^{d}\to R^{d}$ is
${\cal F}\otimes {\cal B}(S^1)\otimes {\cal B}(R^{d})$, ${\cal B}(R^{d})$ measurable and $
g^{\omega}: S^1\times R^{d}\to R^{d}$ is jointly continuous for each $\omega$ and $
g^{\omega}(s,\cdot): R^{d}\to R^{d}$ is differentiable  for each $(\omega, s)\in \Omega\times 
S^1$. To be
 convenient, denote
\begin{eqnarray}\label{(11)}
H^{\omega}(s,y)=(h(s) {\mod} 1,g^{\omega}(s,y)),
\end{eqnarray}
for the skew product map on $S^1\times R^{d}$. We shall
 also define $g^{(n)}$ iteratively by
 $H^{(n),\omega}(s,y)=(h^{(n)}(s) {\mod} 1,g^{(n),\omega}(s,y))$.
It is easy to know that $\hat \theta:
\{\cdots,-3,-2,-1,0,1,2,3,\cdots\} \times \Omega\to \Omega $ is a
$P$-preserving map such that for all $m,n\in
\{\cdots,-3,-2,-1,0,1,2,3,\cdots\} $
\begin{eqnarray}
H^{(m+n),\omega}(s,y)=H^{(m),\hat \theta ^n \omega}(H^{(n),\omega}(s,y)),
\end{eqnarray}
for all $(s,y)\in S^1\times R^{d}$ almost surely.

Let $\pi: R\times R^{d}\rightarrow S^1\times R^{d}$ be the natural covering
$\pi(a,y)= (a\hspace{0.2cm} mod\hspace{0.1cm}1,y)$ and
$\varphi^{\omega}$ a periodic  continuous function. If
$\pi({\rm graph} \hspace{0.2cm}\varphi)$ is invariant under $(H,\hat
\theta)$, that is to say that $H^{\hat \theta ^{-1}\omega}\pi({\rm
graph} \hspace{0.2cm}\varphi^{\hat \theta ^{-1}\omega})=\pi({\rm
graph} \hspace{0.2cm}\varphi^{\omega})$ a.s., it is said that
$\varphi$ is an invariant curve for the skew product (\ref{(11)}).

To analyse the structure of $X^{\omega}$, we pose the following
conditions.

\vskip5pt

{\bf Condition (ii)} {\it Assume that there exists $\delta_1 >0$ and $L_1>0$ such that for any $(s,y)\in X^{\omega}$, there exists a Lipschitz continuous
function $f: S^1\to R^{d}$ with Lipschitz constant $L_1$  and $f(s)=y$ such that  $(s^*,f(s^*))\in X^{\omega}$
when $s^*\in [s-\delta_1,s+\delta_1]$.}
\vskip5pt

Define the $(\delta,\varepsilon)$ neighbourhood of $(s,y)\in X^{\omega}$ and $X^{\omega}$ as
\begin{eqnarray*}
B(s,y,\delta,\varepsilon)=\{(s^{\prime},y^{\prime}): s^{\prime}\in [s-\delta,s+\delta],
||y^{\prime}-f(s^{\prime})||\leq \varepsilon\},
\end{eqnarray*}
and, for any $s\in S^1$,
$${B}_s(\delta,\varepsilon)=\bigcup_{(s,y)\in X^{\omega}}B(s,y,\delta,\varepsilon),$$
and
$${B}(X^{\omega},\varepsilon)=\bigcup_{s\in S^1}B_s(\delta,\varepsilon),$$
respectively. In fact, ${B}_s(\delta,\varepsilon)$ is the $\varepsilon$-neighbourhood
of $X^{\omega}\cap ([s-\delta,s+\delta]\times Y^{\omega})$ in the $y$-direction and ${B}(X^{\omega},\varepsilon)$ is the $\varepsilon$-neighbourhood of $X^{\omega}$ in the $y$-direction.
Note that $\bigcup_{s\in S^1}B_s(\delta,\varepsilon)$ does not depend on $\delta$.
\vskip5pt

{\bf Condition (iii)} {\it  Assume
there exists an $\lambda <1$, $\varepsilon_1 >0$ and an $n_0 \in
 \N$ such that for almost all $\omega \in \Omega$,
\begin{eqnarray}
\label{(13)}
||D_y g^{(n_0),\omega}(s,y)||\leq \lambda {\rm \ for \ all\ } (s,y)\in B(X^{\omega},
\varepsilon_1),
\end{eqnarray}
and
\begin{eqnarray}\label{is4}
c={\rm ess}\sup_{\hskip-15pt
\omega}\sup_{(s,y)\in {B}(X^{\omega}, \varepsilon_1)}||\frac{\partial g^{(n_0),\omega }}{\partial s}(s,y)||<+\infty,
\end{eqnarray}
}

This assumption can be understood as a condition on the amplitude of
Lyapunov exponent of the random map.  It is not
difficult to prove that for all $m\in \N$
\begin{eqnarray}\label{is2}
 H^{mn_0,{\omega}}(B(X^{\omega}, \varepsilon_1))
  \subset
{B}(X^{\hat \theta ^{mn_0}\omega},\varepsilon_1).
 \end{eqnarray}
That is to say that there exists a forward random invariant compact set $B(X^{\omega},\varepsilon_1)$.
By the chain rule, for all $(s,y)\in
B(X^{\omega},\varepsilon_1)$ and $m\in \N$,
\begin{eqnarray}\label{hz11}
||D_y g_{s}^{(mn_0), \omega}(y)||\leq
\lambda^m.
\end{eqnarray}
Moreover, it is easy to see that for any $(s,y)\in X^{\omega}$,
\begin{eqnarray}\label{is3}
 H^{mn_0,\omega}(B(s,y,\delta_1, \varepsilon_1))
  \subset
{B}(H^{mn_0,\omega}(s,y),\delta_1,\varepsilon_1),
 \end{eqnarray}
 and for any $(s^{\prime}, y_1), (s^{\prime}, y_2)\in B(s,y, \delta_1,\varepsilon_1)$,
 \begin{eqnarray*}
 H^{mn_0,\omega}(s^{\prime}, y_1)&=&(s^{\prime}+mn_0\ {\rm mod}\ 1,g_{s^{\prime}}^{(mn_0),\omega}(y_1))\in
 B(H^{mn_0,\omega}(s,y),\delta_1,\varepsilon_1),\\
 H^{mn_0,\omega}(s^{\prime}, y_2)&=&(s^{\prime}+mn_0 \ {\rm mod}\ 1,g_{s^{\prime}}^{(mn_0),\omega}(y_2))\in
 B(H^{mn_0,\omega}(s,y),\delta_1,\varepsilon_1)
 \end{eqnarray*}
 and
 \begin{eqnarray}
 ||g_{s^{\prime}}^{(mn_0),\omega}(y_1)-g_{s^{\prime}}^{(mn_0),\omega}(y_2)||\leq \lambda ^m||y_1-y_2||.
 \end{eqnarray}
 From this, it is easy to see that $X^{\omega}$ actually is a random attractor. The main aim of this section 
 is to prove the following theorem.
 
 \begin{thm}\label{march21a}
Under Conditions (i), (ii), (iii),
the cocycle $\gamma$ has $r^{\omega}$ random periodic solutions of 
period $t_{i}^{*\omega}$ with 
random periodic curve $L_i^{\omega}$
of 
 winding number $\tau_i^{\omega}$, $i=1,2,\cdots,r$. Moreover, for any $t\in R$,
$r^{\theta _{-t}\omega}=r^{\omega}$, and $t_i^{*\theta _{-t}\omega}=t_i^{*\omega}$, $\tau _i^{\theta _{-t}\omega}=\tau _i^{\omega}$ for each $i=1,2,\cdots,r^{\omega}$.
\end{thm}

 We need a series of preparations to prove this theorem.  First, for any $s\in S^{1}$, let's define
$$X_s^\omega=X^\omega \cap (\{s\}\times Y^\omega).$$
Let $G_s^{\omega}$ be a connected component of $X_s^{\omega}$. 
For any $y_1,y_2\in R^{d}$ and $(s,y_1), (s,y_2)\in G_s^{\omega}$, there exist
$y_1^{\prime}, y_2^{\prime}\in R^{d}$ with
$(s^{\prime}, y_1^{\prime}), (s^{\prime}, y_2^{\prime})\in X^{\hat \theta ^{-mn_0}\omega}$ such that
\begin{eqnarray*}
s^{\prime}+mn_0\ {\rm mod} \  1=s, \
g_{s^{\prime}}^{mn_0,\hat \theta ^{-mn_0}\omega}(y_1^{\prime})=y_1,\
g_{s^{\prime}}^{mn_0,\hat \theta ^{-mn_0}\omega}(y_2^{\prime})=y_2,
\end{eqnarray*}
and
\begin{eqnarray}\label{march6a}
||y_1-y_2||\leq \lambda ^m ||y_1^{\prime}-y_2^{\prime}||.
\end{eqnarray}

Let
$|G_s^{\omega}|$ be the largest radius of the set $G_s^{\omega}$. Define
$$
|G^{\omega}|=\sup_{s\in S^1}|G_s^{\omega}|.
$$
Set for $l>0$
$$
\Omega _0^l=\{\omega: |G^{\omega}|\geq l\}.
$$
Then we have
\begin{lem} For any $l>0$,
\begin{eqnarray}
P(\Omega _0^l)=0.
\end{eqnarray}
\end{lem}
{\em Proof.}  Assume the claim of the lemma is not true, i.e. there exists
$\alpha _l>0$ such that $P(\Omega _0^l)=\alpha_l$. Also define
$$
\Omega _m^l=\{\omega: |G^{\omega}|\geq {l\over \lambda ^m}\}.
$$
Then from the invariance of $X^{\omega}$, the Conditions (ii) and (iii) and (\ref{march6a}), we know that
$$
\hat \theta ^{m} \left ((\Omega _m^l)^c\right )\subset (\Omega _0^l)^c.
$$
So
$$
\hat \theta ^{m} \left (\Omega _m^l\right )\supset \Omega _0^l.
$$
Now as $\hat\theta$ is measure preserving we know that for any $m\geq 0$,
\begin{eqnarray}\label{feb22a}
P\left (\Omega _m^l\right )=P\left (\hat \theta ^m\Omega _m^l\right )\geq
P\left (\Omega _0^l\right )=\alpha _l>0.
\end{eqnarray}
Note that $\Omega _m^l\downarrow \emptyset$. So by the continuity of probability measure
with respect to decreasing sequence of sets, we have as $m\to +\infty$
\begin{eqnarray}\label{feb22b}
P\left (\Omega _m^l\right )\to 0.
\end{eqnarray}
Clearly, (\ref{feb22a}) and (\ref{feb22b}) contradict each other. The lemma is proved.
\hfill $\blacksquare$
\medskip

In the following, we will use the pullback of random maps
(\cite{arnold}), the Poincar\'e map, and the Lyapunov exponent
 to prove that  $X^{\omega}$ is a union of finite number
of Lipschitz periodic curves. That is to say that there exists $r^{\omega}$
continuous periodic functions $\varphi_1^{\omega}$, $\varphi_2^{\omega},\cdots,
\varphi_{r^{\omega}}^{\omega}$ on $R^1$ with periods $\tau_1^{\omega},\tau_2^{\omega},\cdots,
\tau_{r^\omega}^{\omega}\in \N$ respectively such that
\begin{eqnarray}
X^{\omega}=L_1^{\omega}\cup
L_2^{\omega}\cup \cdots \cup L_{r^{\omega}}^{\omega},
\end{eqnarray}
 where
 \begin{eqnarray}
 L_i^{\omega}={\rm
graph}(\varphi_i^{\omega})=\{(s\mod 1,\varphi^{\omega}_i(s)): s\in
[0,\tau_i^{\omega})\}, \ i=1,2,\cdots,r^{\omega},
\end{eqnarray}
 are invariant under $(H,
\hat\theta)$.
 Some estimates in the proof of the following lemmas (Lemma
 \ref{l1}-Proposition
 \ref{p1}) are the extension of the results in
 \cite{stark} to the stochastic case. This is not trivial and the pullback
 technique has to be used to make the estimates work.
 Moreover, it is essential to prove that
 \begin{eqnarray}
 \tau_i^{\theta_{-t}\omega}
 =\tau_i^{\omega} \ {\rm and} \  r^{\theta_{-t}\omega}=r^{\omega}
 \end{eqnarray}
  for all $t\in R$. The fact
 that the periodic curves are not the trajectories of the random
 dynamical system, and the fact that the periodic curves are different corresponding to
 different $\omega$, make it difficult to follow the trajectories of the
 random dynamical systems. The essential difficulty in the stochastic case arises
 from the fact that although $X^{\omega}$ is a random invariant set satisfying
 (\ref{is1}), but in general, $X^{\omega}$ and $\gamma ^{\omega}(t) X^{\omega}$
 are different sets. The set $X^{\omega}$ is not invariant in the sense as in the deterministic case.
 This is fundamentally different from the deterministic case. As a consequence,
  the neighbourhood of $X^{\omega}$ is only forward invariant in the sense of
  (\ref{is2}). The assumption (\ref{hz11}) makes the $(\delta_1,\varepsilon_1)$-neighbourhood of $X^{\omega}$ under the map $H^{mn_0, \omega}$ contracting in the $y$-direction to the
  $(\delta_1,\varepsilon_1)$-neighbourhood of $X^{\hat \theta ^{mn_0}\omega}$, rather than
  the $(\delta_1,\varepsilon_1)$-neighbourhood of $X^{\omega}$ as in the deterministic case.
  Locally, (\ref{is3}) says the map $H^{mn_0,\omega}$ maps the
  $(\delta_1,\varepsilon_1)$-neighbourhood of $(s,y)\in X^{\omega}$ to the
  $(\delta_1,\varepsilon_1)$-neighbourhood of $H^{mn_0,\omega}(s,y)\in 
  X^{\hat \theta ^{mn_0}\omega}$. Here, unlike the deterministic case,
  $H^{mn_0,\omega}(s,y)$ is not on the same $X^{\cdot}$ as $(s,y)$ is.

 To prove the above claim, for $(s, y)\in B(X^{\omega},\varepsilon_1)\subset S^1\times Y^{
 \omega},$ denote
\begin{eqnarray*}
{h_1 (s)} &=& {h^{(n_0)} (s)},\\
g_1^{\omega}(s,y)&=& {g^{(n_0),\omega}} (s,y) = g ^{\hat \theta ^{n_0-1}\omega}(h^{(n_0-1)} (s), g ^{(n_0-1),\omega} (s,
y)), \\
H_1^{\omega}(s,y)&=&(h_1(s),g_1^{\omega}(s,y)).
\end{eqnarray*}
For any $(s^*, y^*) \ \in \ S_1\times Y^{\hat
\theta^{-mn_0}\omega}$, define
$$\xi_m^{\hat \theta^{-mn_0}\omega} : [h_1^m (s^*)-\delta_1, h_1^m (s^*) + \delta_1] \to Y^{\omega}$$
by an induction:
\begin{eqnarray*}
\xi_0^{\hat \theta^{-mn_0}\omega}(s)&=&f(s)\in Y^{\hat \theta^{-mn_0}\omega},\forall s\in [s^*-\delta_1,s^*+\delta_1],\\
\\
\xi_1^{\hat \theta^{-mn_0}\omega}(s)&=&g_1^{\hat \theta^{-mn_0}\omega}(h_1^{-1}(s),\ \xi_0^{\hat \theta^{-mn_0}\omega}(h_1^{-1}(s)))\in Y^{\hat \theta ^{-(m-1)n_0}\omega},\\
\\
&&\qquad\qquad \forall  s\in[h_1(s^*)-\delta_1,\
h_1(s^*)+\delta_1],
\end{eqnarray*}
and
\begin{eqnarray*}
\xi_m^{\hat \theta^{-mn_0}\omega}(s)&=&g_1^{\hat \theta^{-n_0}\omega}
(h_1^{-1}(s),\ \xi_{m-1}^{\hat \theta^{-mn_0}\omega}(h_1^{-1}(s)))
\in Y^\omega,\\
\\
&&\quad \forall s\in [h_1^m(s^*)-\delta_1,\
h_1^m(s^*)+\delta_1].
\end{eqnarray*}
Denote
$$\delta_2= min\{\frac{\delta_1}2,\frac{\varepsilon_1}{4L_1}\} \ \ and \ \ L=\frac{c}{1-\lambda}.$$

\begin{lem}\label{l1}
\qquad Under Conditions (i), (ii), (iii), the function $\xi_i^{\hat \theta^{-mn_0}\omega}$ is Lipschitz
continuous
with a Lipschitz constant $L$ for all $m\in \N$ and $i=1,2,\cdots,m$, that is, for any $m\in \N$ and
$i=1,2,\cdots,m$,
\begin{eqnarray*}
(s,\xi_i^{\hat \theta^{-mn_0}\omega}(s)), (s^{\prime},\xi_i^{\hat \theta^{-mn_0}\omega}(s^{\prime}))
\in B(X^{\hat\theta ^{-(m-i)n_0}\omega},\varepsilon_1),
\end{eqnarray*} and
\begin{eqnarray*}
\|\xi_i^{\hat \theta^{-mn_0}\omega}(s)-\xi_i^{\hat \theta^{-mn_0}\omega}(s^{\prime})\|\le L|s-s^{\prime}|,
\end{eqnarray*}
for any $s,\ s^{\prime}\in[h_1^i(s^*)-\delta_2,\ h_1^i(s^*)+\delta_2].$
\end{lem}
{\bf Proof}\qquad We prove this lemma by the induction on $i=1,2,\cdots,m$ for an arbitrary $m$. When $i=1$ and
$s,\ s^{\prime}\in[h_1(s^*)-\delta_2,\ h_1(s^*)+\delta_2]$, by (\ref{is4}), we have
\begin{eqnarray*}
&&
\|\xi_1^{\hat \theta^{-mn_0}\omega}(s)-\xi_1^{\hat \theta^{-mn_0}\omega}(s^{\prime})\|\\
&=& |g_1^{\hat \theta^{-mn_0}\omega}\ (h_1^{-1}(s),\
f(s))-g_1^{\hat \theta^{-mn_0}\omega}(h_1^{-1}(s^{\prime}),\ f(s))\|\\
&\leq & c |h_1^{-1}(s)-h_1^{-1}(s^{\prime})|\\
&<& L|s-s^{\prime}|,
\end{eqnarray*}
and by (\ref{is2}), we have $(s,\xi_1^{\hat
\theta^{-mn_0}\omega}(s)), (s^{\prime},\xi_1^{\hat
\theta^{-mn_0}\omega}(s^{\prime})) \in B(X^{\hat\theta
^{-(m-1)n_0}\omega},\varepsilon_1)$.
 Now suppose the required result holds
for $i-1\in \{1,2,\cdots,m-1\}$, then for any $s,\ s^{\prime}\in[h_1^i(s^*)-\delta_2,\ h_1^i(s^*)+\delta_2]$
\begin{eqnarray*}
&&
\|\xi_i^{\hat \theta^{-mn_0}\omega}(s)-\xi_i^{\hat \theta^{-mn_0}\omega}(s^{\prime})\|\\
&=&\|g_1^{\hat \theta^{-n_0}\omega}\ (h_1^{-1}(s),\ \xi_{i-1}^{\hat
\theta^{-mn_0}\omega}(h_1^{-1}(s)))  -g_1^{\hat \theta^{-n_0}\omega}\
(h_1^{-1}(s^{\prime}),\ \xi_{i-1}^{\hat
\theta^{-mn_0}\omega}(h_1^{-1}(s^{\prime})))||
\\
&\le&|g_1^{\hat \theta^{-n_0}\omega}(h_1^{-1}(s),\ \xi_{i-1}^{\hat
\theta^{-mn_0}\omega}(h_1^{-1}(s)))-g_1^{\hat
\theta^{-n_0}\omega}(h_1^{-1}(s^{\prime}),\
\xi_{i-1}^{\hat \theta^{-mn_0}\omega}(h_1^{-1}(s)))\|\\
&& +\|g_1^{\hat \theta^{-n_0}\omega}(h_1^{-1}(s^{\prime}),\
\xi_{i-1}^{\hat \theta^{-mn_0}\omega}(h_1^{-1}(s))) -g_1^{\hat
\theta^{-n_0}\omega}(h_1^{-1}(s^{\prime}),\ \xi_{i-1}^{\hat
\theta^{-mn_0}\omega}(h_1^{-1}(s^{\prime})))\|
\\
&\le & c |h_1^{-1}(s)-h_1^{-1}(s^{\prime})|
+\lambda \|\xi_{i-1}^{\hat \theta^{-mn_0}\omega}(h_1^{-1}(s))-\xi_{i-1}^{\hat \theta^{-mn_0}\omega}(h_1^{-1}(s^{\prime}))\|\\
&\le& (c+\lambda L)|s-s^{\prime}|\\
&\le & L|s^{\prime}-s|,
\end{eqnarray*}
and the claim $(s,\xi_i^{\hat \theta^{-mn_0}\omega}(s)), (s^{\prime},\xi_i^{\hat \theta^{-mn_0}\omega}(s^{\prime}))
\in B(X^{\hat\theta ^{-(m-i)n_0}\omega}, \varepsilon_1)$ follows from (\ref{is2}).
The lemma is proved.
\hfill $\blacksquare$
\medskip

For any $(s,y)\in X_s^\omega$, let $N(s,y,\delta_2,\varepsilon_1)$
be the interior of $B^\omega(s,y,\delta_2,\varepsilon_1)$.  Then for
any $s^*\in S^1$, $\{N(s^*,y,\delta_2,\varepsilon_1)\big| (s^*,y)\in
X_{s^*}^\omega\}$ is an open covering of $X_{s^*}^\omega$.  By
compactness of $X_{s^*}^\omega$, a finite subcover,
$N(s^*,y^{(1)}_\omega,\delta_2,\varepsilon_1)$,
$N(s^*,y_\omega^{(2)},\ \delta_2,\varepsilon_1)$, $\ldots$, $N(s^*,
y_\omega^{(p_\omega)},\delta_2,\varepsilon_1)$, could be found.
Define
$$N^\omega (s^*,\delta_2,\varepsilon_1)=\bigcup_{i=1}^{p_\omega} N(s^*,y_\omega^{(i)},\delta_2,\varepsilon_1),$$
$$B^\omega (s^*,\delta_2,\varepsilon_1)=\bigcup_{i=1}^{p_\omega} B(s^*,y_\omega^{(i)},\delta_2,\varepsilon_1).$$
Note that $B^\omega(s^*,\delta_2,\varepsilon_1)$ is the closure
of $N^\omega(s^*,\delta_2,\varepsilon_1)$. It is easy to see that $X_{s^*}^{\omega}\in N^{\omega}
(s^*,\delta_2,\varepsilon_1)$.



It is possible for $B(s^*,y^{(i)},\delta_2,\varepsilon_1)$ to
overlap, which leads to the inconvenience in the argument below.  It
is therefore to merge such boxes and work with the connected
components of $B^\omega(s^*,\delta_2,\varepsilon_1)$.  Denote them
by $B_1^\omega(s^*,\delta_2,\varepsilon_1),B_2^\omega(s^*,\delta_2,\varepsilon_1),\ldots ,
B_{r^{*\omega}}^\omega(s^*,\delta_2,\varepsilon_1)$ and let the minimal distance between any
two of them be $\Delta^\omega>0$.  Note that the diameter of any $
B_j^\omega(s^*,\delta_2,\varepsilon_1)$ in the $y$-direction is at most $b^*$.  Later in Lemma \ref{l9} it will be proved that $r^{*\omega}=r^{*\hat
\theta ^{-mn_0}\omega}$. But we don't need this result till the proof of Proposition \ref{p1}.

 \bigbreak

\begin{lem}\label{f6}
\qquad Under Conditions (i), (ii), (iii), for any $j\in\{1,2,\ldots ,r^{*\hat\theta ^{-mn_0}\omega}\}$ and any $m\in \N$,
$$\|y-y^{\prime}\|\le L|s-s^{\prime}|+2\lambda^m b^*,$$
$$\forall (s,y),\ (s^{\prime},y^{\prime})\in H^{m,\hat \theta^{-mn_0}\omega}(B_j^{\hat \theta^{-mn_0}\omega}(s^*,\delta_2,\varepsilon_1)).$$
\end{lem}
{\bf Proof}\qquad Choose $(h_1^{-m}(s), \hat y),\
(h_1^{-m}(s^{\prime}),\hat y^{\prime})\in B_j^{\hat
\theta^{-mn_0}\omega}(s^*,\delta_2,\varepsilon_1)$ such that
\begin{eqnarray*}H_1^{m,\hat \theta^{-mn_0}\omega}(h_1^{-m}(s),\hat y)&=&(s,y),\\
\\
H_1^{m, \hat \theta^{-mn_0}\omega}(h_1^{-m}(s^{\prime}),\hat
y^{\prime})&=&(s^{\prime}, y^{\prime}).\end{eqnarray*}
Then it is obvious that
\begin{eqnarray*}y&=&g_1^{(m),\hat \theta^{-mn_0}\omega}(h_1^{-m}(s),\hat y),\\
y'&=&g_1^{(m), \hat \theta^{-mn_0}\omega}(h_1^{-m}(s^{\prime}),\hat
y^{\prime}).\end{eqnarray*}
 Let $(s^*,y^*)\in
X^{\hat \theta^{-mn_0}\omega}\cap {B}_j^{\hat
\theta^{-mn_0}\omega}(s^*,\delta_2,\varepsilon_1)$ such that $(h_1^{-m}(s),
y^*),(h_1^{-m}(s^{\prime}),y^*) \linebreak
\in {B}_j^{\hat
\theta^{-mn_0}\omega}(s^*,\delta_2,\varepsilon_1)$, then from (\ref{hz11}) and Lemma \ref{l1},
\begin{eqnarray*}
\|y-y'\|&=&\|g_1^{(m),\hat \theta^{-mn_0}\omega}(h_1^{-m}(s),\ \hat
y)  -g_1^{(m),\hat \theta^{-mn_0}\omega}(h_1^{-m}(s^{\prime}),\hat
y^{\prime})\|\\
&\le & \|g_1^{(m),\hat \theta^{-mn_0}\omega}(h_1^{-m}(s),\hat
y)-g_1^{(m),\hat \theta^{-mn_0}\omega}(h_1^{-m}(s),
y^*)\|\\
&& +\|g_1^{(m),\hat \theta^{-mn_0}\omega}(h_1^{-m}(s),
y^*)-g_1^{(m),\hat \theta^{-mn_0}\omega}(h_1^{-m}(s^{\prime}),y^*)\|\\
&& +\|g_1^{(m),\hat
\theta^{-mn_0}\omega}(h_1^{-m}(s^{\prime}),y^*)-g_1^{(m),\hat
\theta^{-mn_0}\omega}(h_1^{-m}(s^{\prime}),\hat
y^{\prime})\|\\
&\le &
 2\lambda^m b^*+L|h_1^{-m}(s)-h_1^{-m}(s^{\prime})|\\
&\le& 2\lambda^m b^*+L|s-s^{\prime}|.
\end{eqnarray*}
 \hfill$\blacksquare$

Choose $N^\omega\in\N$ such that
$$N^\omega>\frac{{\rm log}\,\ds\frac{\Delta^\omega}{2b^*}}{{\rm log}\,\lambda}.$$
This implies that
$$2\lambda^{N^\omega}b^*<\Delta^\omega.$$
Choose $\delta_3\in (0,\delta_2)$ satisfying
\begin{eqnarray}\label{december1}
\delta_3^\omega<\frac{\Delta^\omega-2\lambda^{N^\omega}
b^*}{L}.
\end{eqnarray}

\begin{lem} \label{le3.4}
Under Conditions (i), (ii), (iii), for any $m\in  \{N^\omega, N^\omega+1, N^\omega+2, \dots. \}$
and any  $j\in \left \{1,2, \dots, r^{*\hat \theta
^{-mn_0}\omega}\right \}$,  there exists a unique $i \in  \left \{1,
2,\dots, r^{*\omega}\right \}$ such that
$$
{B_i}^{\omega} (s^*,\delta_3,\varepsilon_1) \cap H_1^{m, \hat
\theta^{-mn_0}\omega} \left ( B_j^{\hat \theta^{-mn_0}\omega
}(s^*,\delta_3,\varepsilon_1) \right ) \neq \emptyset.
$$
\end{lem}
{\bf Proof}\qquad By the definition of ${B}_j^{\hat \theta^{-mn_0}
\omega} (s^*,\delta_3,\varepsilon_1)$, we know that there exists a $ (s^*,
y^*) \in X^{\hat \theta^{-mn_0} \omega} \cap {B}_j^{\hat
\theta^{-mn_0} \omega} (s^*,\delta_3,\varepsilon_1)$. Because of the
invariance of $X$ with respect to $H_1$, we know that $H_1^{m, \hat
\theta^{-mn_0} \omega} (s^*, y^*) \in X^\omega$. 
So $H_1^{m, {\hat \theta^{-mn_0}\omega}}
(s^*, y^*) \in B^\omega (s^*, \delta_3, \varepsilon_1)$. Hence there
exists an $i \in \left \{ 1, 2,\dots, r^{*\omega}\right \}$ such
that $H_1^{m, \hat \theta^{-mn_0} \omega} (s^*, y^*) \in {B}_i^\omega
(s^*,\delta_3,\varepsilon_1).$ So
$$
{B}_i^\omega (s^*,\delta_3,\varepsilon_1) \cap H_1^{m, \hat \theta^{-mn_0}
\omega} \left (B_j ^{\hat \theta^{-mn_0} \omega} (s^*,\delta_3,\varepsilon_1)
\right ) \neq \emptyset.
$$
Now we prove the uniqueness of $i$.  For any $(s, y) \in
 {B}_j^{\hat \theta^{-mn_0}\omega} (s^*,\delta_3,\varepsilon_1)$, $ (h_1^m (s),
g_1^{m,\hat\theta^{-mn_0}\omega} (s, y)) \in H^{m,\hat
\theta^{-mn_0}\omega}
 (B_j^{\hat
\theta^{-mn_0} \omega} (s^*,\delta_3,\varepsilon_1))$.  From Lemma \ref{f6} and (\ref{december1})
we
know that
\begin{eqnarray*}
||g_1^{m,\hat\theta^{-mn_0}\omega} (s^*,
y^*)-g_1^{m,\hat\theta^{-mn_0}\omega}(s, y) || & \leq & L  | s - s^*
| + 2
\lambda ^m b^*\\
& < & L \delta_3^\omega + 2 \lambda ^m  b^*\\
&<& \Delta^\omega .
\end{eqnarray*}
So for any $i^{\prime} \in \{ 1, 2, \dots,r^{*\omega}\}\setminus \{
{i} \}$, $(h_1 ^{(m)} (s), g_1^{m,\hat\theta^{-mn_0}\omega}(s, y))
\notin {B}_{i^{\prime}}^\omega (s^*,\delta_3,\varepsilon_1)$. Thus
$$H_1^{m,\hat \theta^{-mn_0} \omega} ({B}_j^{\hat \theta^{-mn_0} \omega} (s^*,\delta_3,\varepsilon_1))
 \cap B_{i^{\prime}}
^{\omega}(s^*,\delta_3,\varepsilon_1) = \emptyset,
$$
and the uniqueness of $i$ is proved.
\hfill $\blacksquare$
\bigskip

\begin{defi}  Given any $ m\in \{N^\omega, N^\omega+1, N^\omega+2, \dots. \}$ and $ j \in \{1, 2, \dots, r^{*\hat \theta
^{-mn_0}\omega}\}$, denote by $\sigma_m^\omega ({j})$ the unique $i
\in  \{1, 2, \dots, r^{*\omega} \}$ such that \  ${B}_i^{\omega}
(s^*,\delta_3,\varepsilon_1)\cap H_1^{m,\hat \theta^{-mn_0}\omega} \left ({B_j}^{\hat
\theta^{-mn_0} \omega} (s^*,\delta_3,\varepsilon_1) \right ) $
$\neq \emptyset$.
\end{defi}

\begin{lem}\label{l9}
Under Conditions (i), (ii) and (iii), for any $m\in \{N^\omega, N^\omega+1, N^\omega+2, \dots. \}$, $r^{*\hat
\theta ^{-mn_0}\omega}=r^{*\omega}$ and the function
$\sigma_m^{\omega}:\{1,2,\dots,r^{*\omega}\}\to
\{1,2,\dots,r^{*\omega}\}$ is a permutation. In particular,
$\sigma_m^{\omega}$ is invertible and given any $m\in
\{N^\omega, N^\omega+1, N^\omega+2, \dots. \}$, $i\in\{1, 2,\dots,r^{*\omega}\}$, there
exists a unique $j=(\sigma_m^{\omega})^{-1}(i)=\tau_m^{\hat
\theta^{-mn_0}\omega}(i)\in \{1,2,\cdots,r^{*\omega}\}$ such that
$B_i^{\omega}(s^*,\delta_3,\varepsilon_1)\cap H_1^{m,\hat
\theta^{-mn_0}\omega}(B_j^{\hat \theta^{-mn_0}\omega}(
s^*,\delta_3,\varepsilon_1))\neq \emptyset.$
\end{lem}
{\bf Proof}\qquad  Clearly, for any $i\in\{1,2,\dots,r^{*\omega}\}$,
$ B_i^{\omega}(s^*,\delta_3,\varepsilon_1)\cap X_{s^*}^{\omega}\neq\emptyset$. Hence,
using Lemma \ref{le3.4}, we have
\begin{eqnarray*}
B_i^{\omega}(s^*,\delta_3,\varepsilon_1)\cap\Big(\bigcup_{j=1}^{r^{*\hat\theta
^{-mn_0}\omega}} H_1^{m,\hat \theta^{-mn_0}\omega}\big(B_j ^{\hat
\theta^{-mn_0}\omega}(s^*,\delta_3,\varepsilon_1)\big)\Big)\neq\emptyset.
\end{eqnarray*}
Thus the map $\sigma_m:\{1,2,\cdots, r^{*\hat\theta
^{-mn_0}\omega}\}\to \{1,2,\cdots, r^{*\omega}\}$ is onto. We need to prove that
$\sigma_m$ is one-to-one.
As the map $H_1^{m,\hat \theta^{-mn_0}\omega}$ is a contraction in the $y$-direction and
a shift in the $s$ direction, it is evident that for such a $j$ with $\sigma_m^{\omega}(j)=i$,
\begin{eqnarray}
H_1^{m,\hat
\theta^{-mn_0}\omega}(B_j^{\hat \theta^{-mn_0}\omega}(s^*,\delta_3, \varepsilon_1))\subset B_i^{\omega}(s^*,\delta_3, \varepsilon_1).
\end{eqnarray}
But for $j^{\prime}\ne j$, there exists $i^{\prime}\in \{1,2,\cdots, r^{*\omega}\}$ such that
\begin{eqnarray}
H_1^{m,\hat
\theta^{-mn_0}\omega}(B_{j^{\prime}}^{\hat \theta^{-mn_0}\omega}(s^*,\delta
_3, \varepsilon_1))\subset B_{i^{\prime}}^{\omega}(s^*,\delta
_3, \varepsilon_1),
\end{eqnarray}
and
\begin{eqnarray}
B_{i^{\prime}}^{\omega}(s^*,\delta
_3, \varepsilon_1)\cap B_{i}^{\omega}(s^*,\delta
_3, \varepsilon_1)=\emptyset.
\end{eqnarray}
It follows that
\begin{eqnarray}
B_{i}^{\omega}(s^*,\delta
_3, \varepsilon_1)\cap H_1^{m,\hat
\theta^{-mn_0}\omega}(B_{j^{\prime}}^{\hat \theta^{-mn_0}\omega}(s^*,\delta
_3, \varepsilon_1))=\emptyset.
\end{eqnarray}
 Therefore
$\sigma_m$ is an one-to-one map and $r^{*\hat\theta ^{-mn_0}\omega}=r^{*\omega}$. In particular, $\sigma _m^{\omega}$ is a permutation. \hfill$\blacksquare$

\begin{prop}\label{p1} Under Conditions (i), (ii), (iii), there exist $r^{*\omega}$
Lipschitz functions $\varphi _i^\omega : [s^*-\delta_3^\omega,
s^*+\delta_3^\omega]\to Y^\omega$ such that $X^\omega\cap
\Big([s^*-\delta_3^\omega, s^*+\delta_3^\omega]\times
Y^\omega\Big)\subset \bigcup\limits_{i=1}^{r^{*\omega}} {\rm graph}
(\varphi_i^\omega)$ and for each $i\in \{1,2,\cdots,r^{*\omega}\}$, we
have ${\rm graph}(\varphi_i^\omega) \subset
B_i^{\omega}(s^*,\delta_3,\varepsilon_1)$.
\end{prop}
{\bf Proof}: Let $\tau_m ^{\hat \theta^{-mn_0}\omega}$ be the
inverse of $\sigma_m^{\omega}$ and for each $i\in
\{1,2,\cdots,r^{*\omega}\}$, define
\begin{eqnarray}\label{zhaoa1}
W_i^{\omega}=\bigcap_{m = N^\omega}^{\infty}H_1^{m, {\hat
\theta^{-mn_0}\omega}}\Big(B_{\tau_m^ {\hat
\theta^{-mn_0}\omega}(i)}^{{\hat
\theta^{-mn_0}\omega}}(s^*,\delta_3,\varepsilon_1)\Big).
\end{eqnarray}
For any $(s,y)$, $(s',y')\in W_i^\omega$, by Lemma $\ref{f6}$, we
have
\begin{eqnarray*}
||y-y'||\leq L|s-s'|+2\lambda^m b^*,
\end{eqnarray*}
for any $m\in \{N^\omega, N^\omega+1, N^\omega+2, \dots. \}$. Let $m\to \infty$, we get
$||y-y'||\leq L|s-s'|$. That is, each $W_i^\omega$ is contained in
the graph of a Lipschitz function with a Lipschitz constant $L$. Let
\begin{eqnarray}\label{zhaoa2}
{\tilde X_i^\omega}&=&X^\omega \cap B_i^{\omega}(s^*,\delta_3,\varepsilon_1)\nonumber\\
&=&X^\omega \cap\Big([s^*-\delta_3, s^*+\delta_3]\times Y^\omega
\Big)\cap  B_i^{\omega}(s^*,\delta_3,\varepsilon_1).
\end{eqnarray}
It is easy to see that
\begin{eqnarray}\label{zhaoa3}
&&
X^\omega \cap\Big([s^*-\delta_3, s^*+\delta_3]\times Y^\omega
\Big)\nonumber\\
&\subset & \bigcap _{m = N^\omega}^{\infty}\bigcup
_{j=1}^{r^{*\omega}} H_1^{m,{\hat
\theta^{-mn_0}\omega}}\Big(B_j^{{\hat
\theta^{-mn_0}\omega}}(s^*,\delta_3,\varepsilon_1)\Big).
\end{eqnarray}
By Lemma \ref{l9}, for any $m\in \{N^\omega, N^\omega+1, N^\omega+2, \dots. \}$,
\begin{eqnarray}
 B_i^{\omega}(s^*,\delta_3,\varepsilon_1)\cap H^{m,{\hat \theta^{-mn_0}\omega}}\Big(B_{\tau_m^
 {\hat \theta^{-mn_0}\omega}(i)}^{{\hat \theta^{-mn_0}\omega}}(s^*,\delta_3,\varepsilon_1)\Big)\neq\emptyset,
 \label{zhaoa4}\\
 B_i^{\omega}(s^*,\delta_3,\varepsilon_1)\cap\Big(\bigcup_{j\neq \tau_m^ {\hat \theta^{-mn_0}
 \omega}(i)}^{r^{*\omega}}H^{m,{\hat \theta^{-mn_0}\omega}}\Big(B_j^{{\hat \theta^{-mn_0}\omega}}(s^*,\delta_3,\varepsilon_1)\Big)\Big)=\emptyset.\label{zhaoa5}
 \end{eqnarray}
So it follows from (\ref{zhaoa1})-(\ref{zhaoa5}) that
\begin{eqnarray*}
\tilde X_i^{\omega}\subset W_i^\omega.
\end{eqnarray*}
It is easy to show that
\begin{eqnarray*}
\tilde X_i^\omega\cap X_s^\omega\neq \emptyset,
\end{eqnarray*}
for any $s\in [s^*-\delta_3, s^*+\delta_3]$, $i=1,2,\cdots,r^{*\omega}$.
Moreover, for each $s\in[s^*-\delta_3,\ s^*+\delta_3],\ \tilde
{X}_i^\omega\cap X^\omega_s$ contains exactly one point. This can be
seen from $||y-y'||\leq L |s-s'|$, for any $(s,\ y),\ (s',\ y')\in
\tilde {X}^\omega_i.$ Denote this point by $(s,\ \varphi^\omega_i(s)).$
But when $s$ varies in $[s^*-\delta_3,s^*+\delta_3]$,
$(s,\varphi_i^{\omega}(s))$ traces the ${\rm graph}(\varphi_i^{\omega})$.
It is obvious that ${\rm graph}(\varphi^\omega_i)= \tilde
{X}^\omega_i\subset B^\omega_i(s^*,\delta_3, \varepsilon_1),$  and
$\varphi^\omega_i$ is a Lipschitz function with the Lipschitz constant
$L$. \hfill$\blacksquare$
\medskip

\begin{thm}\label{thm1}
Under Conditions (i), (ii), (iii),
$X^\omega$ is a union of a finite number of Lipschitz periodic curves.
\end{thm}
{\bf Proof}: First note that $\delta_3>0$ is independent of $s^*\in
S^1$. Let $M\in \N$ such that ${1\over M}\leq \delta_3$. Define
$s_m={m\over M}$, $m=1,2,\cdots,M$. Then $\{(s_{m-1}, s_{m+1}):
m=1,2,\cdots, M\}$ (in which $s_{M+1}=s_1, s_0=s_M$) covers $S^1$.
By Proposition $\ref{p1}$, we know that $X^\omega\cap
\Big([s_{m-1},s_{m+1}]\times Y^\omega\Big)$ contains a finite number
of Lipschitz curves, denote their number by $r^{*\omega}(m)$. Since
$[s_{m-1},s_m]\subset [s_{m-2},s_m]\cap[s_{m-1},s_{m+1}]$, so we
have $r^{*\omega}(m_1)=r^{*\omega}(m_2)$ when $m_1\neq m_2$. So
$r^{*\omega}$ is independent of $m$ and define all of them by
$r^{*\omega}$. Thus the Lipschitz curves on $X^\omega\cap
\Big([s_{m-1},s_{m+1}]\times Y^\omega\Big)$ could be expanded to
$S^1$ and we have the following random Poincare map
$$
H^{n_0,\hat\theta ^{-n_0}\omega}:G_{s}^{\hat\theta ^{-n_0}\omega}\to G_s^{\omega},
$$
in which $G_{s}^{\hat\theta ^{-n_0}\omega}, G_s^{\omega}$ are finite sets containing $r^{*\omega}$ elements:
\begin{eqnarray*}
G_{s}^{\hat\theta ^{-n_0}\omega}&=&\{(s \mod 1, \varphi_i^{\hat \theta ^{-n_0}\omega}(s)): i=1,2,\cdots,r^{*\omega}\},\\
G_s^\omega&=&\{(s\mod 1, \varphi_i^{\omega}(s)): i=1,2,\cdots,r^{*\omega}\},
\end{eqnarray*}
for a fixed $s\in R^1$.
By the finiteness of $G_s^{\omega}$, we know
\begin{eqnarray*}
\varphi_i^{\omega}(s+1)&=&\varphi_{i_1}^{\omega}(s),\\
\varphi_i^{\omega}(s+2)&=&\varphi_{i_2}^{\omega}(s),\\
&&\cdots\\
\varphi_i^{\omega}(s+r^{*\omega}))
&=&\varphi_{i_{r^{*\omega}}}^{\omega}(s).
\end{eqnarray*}
Actually above is true for any $s$ due to the continuity of
$\varphi_i^{\omega}$, $i=1,2,\cdots,r^{*\omega}$. Therefore there are three cases:-\\

(i). Exact one of $i_1,i_2,\cdots, i_{r^{*\omega}}$ is equal to $i$.
Say $i_{\tau_i^{\omega}}=i$. Then
\begin{eqnarray*}
\varphi_i^{\omega}(s+\tau_i^{\omega})=\varphi_i^{\omega}(s),
\end{eqnarray*}
for any $s\in R$.
So $\varphi_i^{\omega}$ is a periodic function of period $\tau_i^{\omega}$.\\

(ii). More than one of $i_1,i_2,\cdots ,i_{r^{*\omega}}$ is equal to
$i$. Denote $\tau_i^{\omega}$ the smallest number $j$ such that
$i_j=i$ and $\tilde \tau_i^{\omega}>\tau_i^{\omega}$ such that
$i_{\tilde\tau_i^{\omega}}=i$. Then
\begin{eqnarray*}
\varphi_i^{\omega}(s+\tau_i^{\omega})&=&\varphi_i^{\omega}(s),\\
\varphi_i^{\omega}(s+\tilde\tau_i^{\omega})&=&\varphi_i^{\omega}(s).
\end{eqnarray*}
But
\begin{eqnarray*}
\varphi_i^{\omega}(s+\tilde\tau_i^{\omega})&=&\varphi_i^{\omega}(s+\tilde{\tau}_i^{\omega}-\tau_i^{\omega}+\tau_i^{\omega})\\
&=&\varphi_i^{\omega}(s+\tilde{\tau}_i^{\omega}-\tau_i^{\omega})\\
&=&\cdots\\
&=&\varphi_i^{\omega}(s+\tilde{\tau}_i^{\omega}-k\tau_i^{\omega}),
\end{eqnarray*}
where $k$ is the smallest integer such that $\tilde
\tau_i^{\omega}-(k+1)\tau_i^{\omega} \leq 0$. Then by definition of $\tau_i^{\omega}$,
$$\tilde \tau_i^{\omega}-{k}\tau_i^{\omega}= \tau_i^{\omega},$$ so
$$\tilde \tau_i^{\omega}={(k+1)}\tau_i^{\omega}.
$$
Therefore $\varphi_i$ is a periodic function of period $\tau_i^{\omega}$.\\

(iii). None of $i_1,i_2,\cdots,i_{r^{*\omega}}$ is equal to $i$. In this case,
at least two of $i_1,i_2,\cdots,i_{r^{*\omega}}$ must be equal. Say
$\tau_2^{\omega}>\tau_1^{\omega}$ are the two such integers such that
$i_{\tau_1^{\omega}}=i_{\tau_2^{\omega}}$ with smallest difference
$\tau_2^{\omega}-\tau_1^{\omega}$. Then
$$\varphi_{i}^{\omega}(s+\tau_1^{\omega})=\varphi_{i}^{\omega}(s+\tau_2^{\omega}).$$
Denote $s+\tau_1^{\omega}$ by $s_{1}$, then
$$\varphi_{i}^{\omega}(s)=\varphi_{i}^{\omega}(s+\tau_2^{\omega}-\tau_1^{\omega}),\forall s\in\mathbb{R}^{1}.$$
Same as (ii) we can see for all other possible $\tilde{\tau}_{2}^{\omega}$
and $\tilde{\tau}_{1}^{\omega}$, $\tilde{\tau}_{2}^{\omega}>\tilde{\tau}_{1}^{\omega}$ and
$i_{\tilde{\tau}_{2}^{\omega}}=i_{\tilde{\tau}_{1}^{\omega}}$,
$\tilde{\tau}_{2}^{\omega}-\tilde{\tau}_{1}^{\omega}$ must be an integer multiple of
$\tau_2^{\omega}-\tau_1^{\omega}$. That is to say $\varphi_{i}^{\omega}$ is a periodic curve
with period $\tau_{2}^{\omega}-\tau_{1}^{\omega}$. \hfill
$\blacksquare$
\medskip

Theorem \ref{thm1} says there exists a finite number of continuous periodic
functions $\varphi_{1}^{\omega},\varphi_{2}^{\omega},\cdots,\varphi_{r^{\omega}}^{\omega}$ on $\mathbb{R}^{1}$. Denote their periods by $\tau_{1}^{\omega},\tau_{2}^{\omega},\cdots,\tau_{r^{\omega}}^{\omega}\in\mathbb{N}$
respectively. So
$$X^{\omega}=L^{\omega}_{1}\cup L^{\omega}_{2}\cup\cdots\cup L^{\omega}_{r^{\omega}}$$
where
$$L^{\omega}_{i}={\rm graph}(\varphi^{\omega}_{i})=\{(s\ mod \ 1, \varphi^{\omega}_{i}(s)): \ s\in[0,\tau_{i}^{\omega})\}.$$
But
\begin{eqnarray*}
H_1^{\hat \theta ^{-n_0}\omega}(X^{\hat \theta ^{-n_0}\omega})=X^{\omega}.
\end{eqnarray*}
So
\begin{eqnarray}\label{is5}
H_1^{\hat \theta ^{-n_0}\omega}(L^{\hat \theta ^{-n_0}\omega}_{1})
\cup H_1^{\hat \theta ^{-n_0}\omega}(L^{\hat \theta ^{-n_0}\omega}_{2})\cup \cdots\cup H_1^{\hat \theta ^{-n_0}\omega}
(L^{\hat \theta ^{-n_0}\omega}_{r^{\hat \theta ^{-n_0}\omega}}) =L^{\omega}_{1}\cup L^{\omega}_{2}\cup\cdots\cup L^{\omega}_{r^{\omega}}.
\end{eqnarray}
It is easy to know that $H_1^{\hat \theta ^{-n_0}\omega}(L^{\hat \theta ^{-n_0}\omega}_i)$ is a closed curve since $L^{\hat \theta ^{-n_0}\omega}_i$ is a closed curve and $H_1^{\hat \theta ^{-n_0}{\omega}}$ is a continuous map. Moreover, since $H_1$ is a homeomorphism, so
\begin{eqnarray}
H_1^{\hat \theta ^{-n_0}\omega}(L^{\hat \theta ^{-n_0}\omega}_i)\cap
H_1^{\hat \theta ^{-n_0}\omega}(L^{\hat \theta ^{-n_0}\omega}_j)=\emptyset,
\end{eqnarray}
when $i\ne j$. Therefore the left hand side of of (\ref{is5}) is a union of $r^{\hat\theta^{-n_0}\omega}$ distinct closed curves and the right hand side of of (\ref{is5}) is a union of $r^{\omega}$
distinct closed curves. Thus for any $i\in\{1,2,\cdots,r^{\hat\theta^{-n_0}\omega}\}$, there exists a unique
$j\in\{1,2,\cdots,r^{\omega}\}$ such that
\begin{eqnarray}
H_1^{\hat \theta ^{-n_0}\omega}(L^{\hat \theta ^{-n_0}\omega}_i)=L_j^{\omega}.
\end{eqnarray}
Denote $j=K(i)$. It is easy to see now that $r^{\hat\theta ^{-n_0}\omega}=r^{\omega}$.
Therefore $K: \{1,2,\cdots,r^{\hat\theta_{-n_0}\omega}\}\to \{1,2,\cdots,r^{\hat\theta_{-n_0}\omega}\}$
and $K$ is a permutation. Reorder $\{L_i^{\omega}: i=1,2,\cdots,r^{\hat\theta_{-n_0}\omega}\}$, we can have
\begin{eqnarray}
H_1^{\hat \theta ^{-n_0}\omega}(L^{\hat \theta ^{-n_0}\omega}_i)=L_i^{\omega}.
\end{eqnarray}
Moreover, using a similar argument, and note for any $t\in R$
$$X^{\theta _{-t}\omega}=L^{\theta _{-t}\omega}_{1}\cup L^{\theta _{-t}\omega}_{2}\cup\cdots\cup L^{\theta _{-t}\omega}_{r^{\theta _{-t}\omega}}$$
and
\begin{eqnarray*}
\gamma^{\theta _{-t}\omega}(t)(L^{\theta _{-t}\omega}_{1})
\cup \gamma^{\theta _{-t}\omega}(t)(L^{\theta _{-t}\omega}_{2})\cup \cdots\cup \gamma^{\theta _{-t}\omega}(t)
(L^{\theta _{-t}\omega}_{r^{\theta _{-t}\omega}}) =L^{\omega}_{1}\cup L^{\omega}_{2}\cup\cdots\cup L^{\omega}_{r^{\omega}}.
\end{eqnarray*}
 we have the following proposition.

\begin{prop}\label{zhaoh1}
Under Conditions (i), (ii) and (iii), for each $\omega$, we have for any $t\in R$
\begin{eqnarray}\label{april28a}
\gamma^{\theta_{-t}\omega}(t)L_i^{\theta _{-t}\omega} = L^{\omega}_{i}.
\end{eqnarray}
\end{prop}

\begin{lem} \label{revision3}
For any $t\in R$, $\tau _i^{\theta _{-t}\omega}=\tau
_i^{\omega}$ for any $i=1,2,\cdots,r^{\omega}$.
\end{lem}
{\bf Proof}. Consider first the case when $t=kt_1 \  \  (k \in  \{0,1,2,\cdots. \})$. Note for any  $s \in  \{0,1,2,\cdots, \tau _i^{\theta _{-t}\omega}\}$,
\begin{eqnarray}\label{april28b}
s_1=\left (\gamma ^{\theta _{-t}\omega}_{(s \mod 1, \varphi_i^{\theta _{-t}\omega}(s))}(t)
\right )_{S_1}
=s+k.
\end{eqnarray}
Here  $(\cdot)_{S_1}$ denotes the $S_1$ coordinate of the vector. So for $t=kt_1$, from (\ref{april28a})
and (\ref{april28b}), it turns out that 
\begin{eqnarray}
\tau _i^{\omega}&=&\left (\gamma ^{\theta _{-t}\omega}_{(\tau _i^{\theta _{-t}\omega} \mod 1, \varphi_i^{\theta _{-t}\omega}(\tau _i^{\theta _{-t}\omega}))}(t)\right )_{S_1}-
\left (\gamma ^{\theta _{-t}\omega}_{(0, \varphi_i^{\theta _{-t}\omega}(0))}(t)\right )_{S_1}\nonumber\\
&=&
\tau
_i^{\theta _{-t}\omega}+k-(0+k)\nonumber\\
&=&
\tau
_i^{\theta _{-t}\omega}.
\end{eqnarray}

Now we consider the case when  $t \in (kt_1, (k+1)t_1)\  \  (k \in  \{0,1,2,\cdots\})$. Note for any  $s \in  \{0,1,2,\cdots, \tau _i^{\theta _{-t}\omega}\}$,
\begin{eqnarray}
s_1=\left (\gamma ^{\theta _{-t}\omega}_{(s \mod 1, \varphi_i^{\theta _{-t}\omega}(s))}(t)\right )_{S_1} \in  (k+s, s+k+1),
\end{eqnarray}
 and for each $(t,\omega)\in (-\infty,+\infty)\times \Omega$, the mapping
$\gamma _{\cdot}^{\omega}(t): S^1\times R^{d}\to S^1\times
R^{d}$ is a homeomorphism. By the periodicity of $\varphi_i^{\theta_{-t}\omega}$ (period being $\tau_i^{\theta _{-t}\omega}$)
\begin{eqnarray}
\left (\gamma ^{\theta _{-t}\omega}_{(\tau _i^{\theta _{-t}\omega} \mod 1, \varphi_i^{\theta _{-t}\omega}(\tau _i^{\theta _{-t}\omega}))}(t)\right )_{S_1}
=\tau _i^{\theta _{-t}\omega}+
\left (\gamma ^{\theta _{-t}\omega}_{0, \varphi_i^{\theta _{-t}\omega}(0)}(t)\right )_{S_1}.
\end{eqnarray}
Now from Proposition 3.9,
\begin{eqnarray}
\tau _i^{\omega}=
\left (\gamma ^{\theta _{-t}\omega}_{(\tau _i^{\theta _{-t}\omega} \mod 1, \varphi_i^{\theta _{-t}\omega}(\tau _i^{\theta _{-t}\omega}))}(t)\right )_{S_1}-
\left (\gamma ^{\theta _{-t}\omega}_{0, \varphi_i^{\theta _{-t}\omega}(0)}(t)\right )_{S_1}=
\tau_i^{\theta _{-t} \omega}.
\end{eqnarray}
 The lemma is proved.
 \hfill$\blacksquare$
 \medskip

{\bf Proof of Theorem \ref{march21a}}.\quad
 First, from (\ref{april28a}), we know that
 once $L^{\omega}_{i}$ is known for one
$\omega$, then $L^{\theta_{-t}\omega}_{i}$is determined for any
$t\in R$.  In the following $\varphi$ is used to represent any $\varphi_{i}$
and $\tau^{\omega}$ the corresponding $\tau_{i}^{\omega}$.  From Lemma \ref{revision3}, 
we know $\tau ^{\theta _{-t}\omega}=\tau
^{\omega}$ for any $t\in R$. Define
\begin{eqnarray}
t^{*}=\tau t_1.
\end{eqnarray}
Then for any $s\in [0,\tau^{\omega})$
\begin{eqnarray}\label{zhaoh2}
\gamma^{\theta_{-t^{*}}\omega}_{(s\ {\rm mod} \ 1,\varphi^{\theta_{-t^{*}}\omega}(s))}(t^{*})=(s\ {\rm mod} \ 1, \varphi^{\omega}(s)).
\end{eqnarray}
Therefore from (\ref{zhaoh2}) and the
cocyle property of $\gamma$, it follows that for any $s\in [0,\tau^{\omega})$
\begin{eqnarray}
\gamma_{(s\ {\rm mod} \
1,\varphi^{\theta_{-t^*-t}\omega}(s))}^{\theta_{-t^*-t}\omega}(t+t^*)\nonumber
&=& \gamma_{\gamma ^{\theta _{-t^*-t}\omega}_{(s\ {\rm mod} \
1,\varphi^{\theta_{-t^*-t}\omega}(s))}(t^*)}^{\theta_{-t}\omega}(t)\nonumber\\&=&
\gamma_{(s\ {\rm mod} \
1,\varphi^{\theta_{-t}\omega}(s))}^{\theta_{-t}\omega}(t).
\end{eqnarray}
This gives that for any $s\in [0,\tau^{\omega})$
\begin{eqnarray}
\gamma_{(s\ {\rm mod} \
1,\varphi^{\omega}(s))}^{\omega}(t+t^*)=\gamma_{(s\
{\rm mod} \ 1,\varphi^{\theta _{t^*}\omega}(s))}^{\theta _{t^*}\omega}(t) \end{eqnarray} for any
$t\leq 0$. That is to say $\gamma$ has a periodic curve
with period $t^*$ and winding number $\tau$. There are $r$ such $\varphi$. That is to say
$\gamma$ has $r$ random periodic solutions. \hfill$\blacksquare$

\begin{acknowledgement}

The authors are very grateful to A. Truman
   who went through the paper and gave many valuable suggestions.
 ZZ acknowledges financial supports of the NSF of China (No. 10371123)
   and National 973 Project (2005 CB 321902), especially
   of  the Royal Society London and the Chinese Academy of Sciences
   that enabled him to visit Loughborough University where the research was done.
   HZ  acknowledges the EPSRC (UK) grant GR/R69518.
   We would like to thank the referee for very useful comments.

\end{acknowledgement}

\end{document}